\documentclass{article}
\usepackage{amsmath}
\usepackage{amssymb}

\usepackage{theorem}
\numberwithin{equation}{section}

\newtheorem{theorem}{Theorem}[section]

\newtheorem{proposition}[theorem]{Proposition}

\begin{document}
\title{Global well-posedness for the Schr{\"o}dinger map problem with small Besov norm}
\date{\today}
\author{Benjamin Dodson}
\maketitle

\noindent \textbf{Abstract:} In this paper we prove a global result for the Schr{\"o}dinger map problem with initial data with small Besov norm at critical regularity.

\section{Introduction}
In this paper we consider the Schr{\"o}dinger map initial value problem

\begin{equation}\label{1.1}
\aligned
\partial_{t} \phi &= \phi \times \Delta \phi, \hspace{5mm} \text{on} \hspace{5mm} \mathbf{R} \times \mathbf{R}^{d}, \\
\phi(0) &= \phi_{0},
\endaligned
\end{equation}

\noindent in dimensions $d \geq 3$, where $\phi : \mathbf{R}^{d} \times \mathbf{R} \rightarrow S^{2} \hookrightarrow \mathbf{R}^{3}$ is a continuous function. The Schr{\"o}dinger map equation arises in ferromagnetism as the Heisenberg model for the ferromagnetic spin system.\vspace{5mm}

\noindent In this paper we are concerned with the issue of global well - posedness of $(\ref{1.1})$ for data which is small in a critical space under the scaling.\vspace{5mm}

\noindent Observe that for any $\lambda > 0$, if $\phi(t, x)$ solves $(\ref{1.1})$, then

\begin{equation}\label{1.2}
\phi(\lambda^{2} t, \lambda x)
\end{equation}

\noindent also solves $(\ref{1.1})$. Now for $Q \in S^{2}$, $\mathbf{Z}_{+} = \{ 0, 1, 2, ... \}$, and $\sigma \in [0, \infty)$, define the space

\begin{equation}\label{1.3}
H_{Q}^{\sigma} = \{ f : \mathbf{R}^{d} \rightarrow \mathbf{R}^{3}, \hspace{5mm} |f(x)| \equiv 1 \hspace{5mm} a.e. \hspace{5mm} f - Q \in H^{\sigma} \},
\end{equation}

\noindent where $H^{\sigma}$ is the usual inhomogeneous Sobolev space. This metric has the induced distance

\begin{equation}\label{1.4}
d_{Q}^{\sigma}(f, g) = \| f - g \|_{H^{\sigma}(\mathbf{R}^{d})}.
\end{equation}

\noindent Also define the metric spaces

\begin{equation}\label{1.5}
H^{\infty} = \cap_{\sigma \in \mathbf{Z}_{+}} H^{\sigma}, \hspace{5mm} \text{and} \hspace{5mm} H_{Q}^{\infty} = \cap_{\sigma \in \mathbf{Z}_{+}} H_{Q}^{\sigma}.
\end{equation}

\noindent The scaling law $(\ref{1.2})$ preserves the $\dot{H}^{d/2}$ and $\dot{H}_{Q}^{d/2}$ homogeneous Sobolev norms, where

\begin{equation}\label{1.6}
\| f \|_{\dot{H}^{\sigma}} = \| \mathcal F(f)(\xi) \cdot |\xi|^{\sigma} \|_{L^{2}},
\end{equation}

\begin{equation}\label{1.7}
\| f \|_{\dot{H}_{Q}^{\sigma}} = \| f - Q \|_{\dot{H}^{\sigma}},
\end{equation}

\noindent and $\mathcal F : L^{2} \rightarrow L^{2}$ is the usual Fourier transform

\begin{equation}\label{1.8}
\mathcal F(f)(\xi) = c_{d} \int_{\mathbf{R}^{d}} e^{-ix \cdot \xi} f(x) dx.
\end{equation}

\noindent $(\ref{1.1})$ is globally well - posed for data which is small in $\dot{H}^{d/2}$.

\begin{theorem}[Global regularity]\label{t1.1}
Assume $d \geq 2$ and $Q \in S^{2}$. Then there exists $\epsilon_{0}(d) > 0$ such that for any $\phi_{0} \in H_{Q}^{\infty}$ with $\| \phi_{0} - Q \|_{\dot{H}^{d/2}} \leq \epsilon_{0}(d)$, then there is a unique solution

\begin{equation}\label{1.9}
\phi = S_{Q}(\phi_{0}) \in C(\mathbf{R} : H_{Q}^{\infty})
\end{equation}

\noindent of the initial value problem $(\ref{1.1})$. Moreover,

\begin{equation}\label{1.10}
\sup_{t \in \mathbf{R}} \| \phi(t) - Q \|_{\dot{H}^{d/2}} \leq C \| \phi_{0} - Q \|_{\dot{H}^{d/2}},
\end{equation}

\noindent and for any $T \in [0, \infty)$ and $\sigma \in \mathbf{Z}_{+}$,

\begin{equation}\label{1.11}
\sup_{t \in [-T, T]} \| \phi(t) \|_{H_{Q}^{\sigma}} \leq C(\sigma, T, \| \phi \|_{H_{Q}^{\sigma}}).
\end{equation}
\end{theorem}

\noindent This theorem was proved in dimensions $d \geq 4$ in \cite{BIK}, and then for dimensions $d \geq 2$ in \cite{BIKT}. \cite{BIKT} also proved a uniform global bound.

\begin{theorem}[Uniform bounds and well - posedness]\label{t1.2}
Assume $d \geq 2$, $Q \in S^{2}$, and $\sigma_{1} \geq \frac{d}{2}$. Then there exists $0 < \epsilon_{0}(d, \sigma_{1}) \leq \epsilon_{0}(d)$ such that if $\phi \in H_{Q}^{\infty}$ with $\| \phi - Q \|_{\dot{H}^{d/2}} \leq \epsilon_{0}(d, \sigma_{1})$, then the global solution $\phi(t)$ to $(\ref{1.1})$ with initial data $\phi_{0}$ satisfies

\begin{equation}\label{1.12}
\sup_{t \in \mathbf{R}} \| \phi(t) - Q \|_{H^{\sigma}} \leq C_{\sigma} \| \phi_{0} - Q \|_{H^{\sigma}}, \hspace{5mm} \frac{d}{2} \leq \sigma \leq \sigma_{1}.
\end{equation}

\noindent In addition, the solution operator admits a continuous extension from

\begin{equation}\label{1.13}
B_{\epsilon_{0}}^{\sigma} = \{ \phi \in \dot{H}_{Q}^{d/2 - 1} \cap \dot{H}^{\sigma} : \| \phi - Q \|_{\dot{H}^{d/2}} \leq \epsilon_{0} \}
\end{equation}

\noindent to $C(\mathbf{R} ; \dot{H}^{\sigma} \cap \dot{H}_{Q}^{d/2 - 1})$.
\end{theorem}

\noindent The case when $d = 2$ is particularly interesting, due to the fact that a solution to $(\ref{1.1})$ conserves the quantities

\begin{equation}\label{1.14}
E_{0}(t) = \int_{\mathbf{R}^{d}} |\phi(t) - Q|^{2} dx,
\end{equation}

\noindent and

\begin{equation}\label{1.15}
E_{1}(t) = \int_{\mathbf{R}^{d}} \sum_{m = 1}^{d} |\partial_{m} \phi(t,x)|^{2} dx.
\end{equation}

\noindent However, in general the proofs in theorems $\ref{t1.1}$ and $\ref{t1.2}$ are more difficult in lower dimensions. Indeed, \cite{BIKT} utilized the caloric gauge of \cite{Taowm} to analyze dimension $d = 2$, because the Coulomb gauge used in \cite{BIKT} was not strong enough.\vspace{5mm}

\noindent In this paper we will extend theorems $\ref{t1.1}$ and $\ref{t1.2}$ to data which is small in Besov - type norms.

\begin{theorem}\label{t1.3}
Suppose $\phi_{0} \in \dot{H}^{d/2}$ and suppose that $\psi(x) \in C_{0}^{\infty}(\mathbf{R}^{d})$ is a radial function supported on $\frac{1}{4} \leq |x| \leq 4$ and $\phi(x) = 1$ on $\frac{1}{2} \leq |x| \leq 2$. Furthermore, suppose that for some $\epsilon_{0}(d, \| \phi_{0} \|_{\dot{H}^{d/2}})$,

\begin{equation}\label{1.16}
\sup_{j} \| \psi(2^{j} \xi) |\xi|^{d/2} \cdot \mathcal F(\phi - Q)(\xi)  \|_{L^{2}(\mathbf{R}^{d})} \leq \epsilon_{0}(d, \| \phi_{0} \|_{\dot{H}^{d/2}}).
\end{equation}

\noindent Then the results of theorems $\ref{t1.1}$ and $\ref{t1.2}$ hold.
\end{theorem}

\noindent In this paper we follow \cite{B}, \cite{BIK}, \cite{BIKT}, \cite{Krieger}, \cite{Krieger1}, \cite{NSU}, and \cite{ShSt} and prove a priori bounds on the derivative of the Schr{\"o}dinger map rather than the Schr{\"o}dinger map itself. We use the Coulomb gauge, as is typically used in analysis of higher dimensions, see for example \cite{BIK}.\vspace{5mm}

\noindent The main new ingredient to this paper is the bilinear estimates obtained from the interaction Morawetz estimates. In \cite{Dodson}, applying the bilinear estimates of \cite{PV} to the various Schr{\"o}dinger map gauges gave bilinear estimates for solutions to $(\ref{1.1})$ in dimensions $d \geq 2$. These results were used in \cite{DS} to improve the results of \cite{Smith1}, \cite{Smith2}, and \cite{Smith3}.\vspace{5mm}

\noindent Here we are able to use bilinear arguments in order to prove theorem $\ref{t1.3}$. These arguments greatly simplify the function spaces that are needed in the analysis.

\section{Bilinear Virial estimates}
The work in this paper will heavily rely on bilinear estimates for solutions to the linear Schr{\"o}dinger equation

\begin{equation}\label{2.1}
i u_{t} + \Delta u = 0.
\end{equation}

\begin{theorem}[Interaction Morawetz estimate]\label{t2.1}
If $u$ solves $(\ref{2.1})$ then

\begin{equation}\label{2.2}
\| |\nabla|^{\frac{3 - d}{2}} |u|^{2} \|_{L_{t,x}^{2}(\mathbf{R} \times \mathbf{R}^{d})}^{2} \lesssim \| u \|_{L_{t}^{\infty} L_{x}^{2}(\mathbf{R} \times \mathbf{R}^{d})}^{2} \| u \|_{L_{t}^{\infty} \dot{H}_{x}^{1/2}(\mathbf{R} \times \mathbf{R}^{d})}^{2}.
\end{equation}
\end{theorem}

\noindent \emph{Proof:} This theorem was proved in \cite{CKSTT2} when $d = 3$. \cite{TVZ} subsequently extended the result to dimensions $d \geq 4$. \cite{PV} and \cite{CGT} independently proved theorem $\ref{t2.1}$ in dimensions $d = 1, 2$. $\Box$\vspace{5mm}

\noindent Of course, \cite{CKSTT2}, \cite{TVZ}, \cite{PV}, and \cite{CGT} also proved theorem $\ref{t2.1}$ for nonlinear Schr{\"o}dinger equations as well. \cite{PV} also made some very important insights into the behavior of linear solutions as well. In particular, \cite{PV} observed that one may compute the interaction Morawetz estimate for two different solutions to $(\ref{2.1})$, yielding a bilinear estimate.

\begin{theorem}[Bilinear virial estimate]\label{t2.2}
Suppose $u$ and $v$ are solutions to $(\ref{2.1})$. Then

\begin{equation}\label{2.3}
\aligned
\int \int_{x_{1} = y_{1}} |\partial_{1}(u(t, x_{1}, ..., x_{d}) \overline{v(t, x_{1}, y_{2}, ..., y_{d})})|^{2} dx dy dt \\ \lesssim \| u \|_{L_{t}^{\infty} \dot{H}^{1/2}(\mathbf{R} \times \mathbf{R}^{d})}^{2} \| v \|_{L_{t}^{\infty} L_{x}^{2}(\mathbf{R} \times \mathbf{R}^{d})}^{2} +  \| v \|_{L_{t}^{\infty} \dot{H}^{1/2}(\mathbf{R} \times \mathbf{R}^{d})}^{2} \| u \|_{L_{t}^{\infty} L_{x}^{2}(\mathbf{R} \times \mathbf{R}^{d})}^{2}.
\endaligned
\end{equation}
\end{theorem}

\noindent In particular, suppose that $u$ and $v$ are solutions to $(\ref{2.1})$ and $P_{j}$ is the usual Littlewood - Paley projection operator. That is, let $\psi(x) = 1$ for $|x| \leq 1$ and $\psi(x) = 0$ for $|x| > 2$ be a smooth, radially symmetric function, decreasing as $r \rightarrow \infty$. Then let

\begin{equation}\label{2.4}
\psi_{j}(x) = \psi(2^{-j} x) - \psi(2^{-j + 1} x),
\end{equation}

\noindent and let $P_{j}$ be the Fourier multiplier given by

\begin{equation}\label{2.5}
\mathcal F(P_{j} f)(\xi) = [\psi(2^{-j} \xi) - \psi(2^{-j + 1} \xi)] \mathcal F(f)(\xi).
\end{equation}

\noindent Then by elementary Fourier analysis,

\begin{equation}\label{2.6}
(P_{j} f)(x) = 2^{jd} \int \check{\psi}(2^{j}(x - y)) f(y) dy,
\end{equation}

\noindent where $\check{\psi}$ is a smooth function that is rapidly decreasing for $|x|$ large, that is, for any integer $N$,

\begin{equation}\label{2.7}
|\check{\psi}(x)| \lesssim_{N} (1 + |x|)^{-N}.
\end{equation}

\noindent Since the Littlewood - Paley operator is a Fourier multiplier, if $u$ solves $(\ref{2.1})$ then $P_{j} u$ also solves $(\ref{2.1})$, so

\begin{equation}\label{2.8}
\aligned
\int \int_{x_{1} = y_{1}} |\partial_{1}(P_{j} u(t, x_{1}, x_{2}, ..., x_{d}) \overline{P_{k} v(t, x_{1}, y_{2}, ..., y_{d})})|^{2} dx dy dt \\ \lesssim \| P_{j} u \|_{L_{t}^{\infty} \dot{H}^{1/2}(\mathbf{R} \times \mathbf{R}^{d})}^{2} \| P_{k} v \|_{L_{t}^{\infty} L_{x}^{2}(\mathbf{R} \times \mathbf{R}^{d})}^{2} +  \| P_{k} v \|_{L_{t}^{\infty} \dot{H}^{1/2}(\mathbf{R} \times \mathbf{R}^{d})}^{2} \| P_{j} u \|_{L_{t}^{\infty} L_{x}^{2}(\mathbf{R} \times \mathbf{R}^{d})}^{2}.
\endaligned
\end{equation}

\noindent But then by $(\ref{2.6})$, $(\ref{2.7})$, $(\ref{2.8})$, and H{\"o}lder's inequality, if $k \leq j$,

\begin{equation}\label{2.9}
\| \partial_{1} (P_{j} u \overline{P_{k} v}) \|_{L_{t,x}^{2}(\mathbf{R} \times \mathbf{R}^{d})}^{2} \lesssim 2^{(d - 1)k} 2^{j} \| P_{j} u \|_{L_{t}^{\infty} L_{x}^{2}(\mathbf{R} \times \mathbf{R}^{d})}^{2} \| P_{k} v \|_{L_{t}^{\infty} L_{x}^{2}(\mathbf{R} \times \mathbf{R}^{d})}^{2}.
\end{equation}

\noindent Now since there is absolutely nothing special about the direction $e_{1} \in S^{d - 1}$, averaging over all directions $\omega \in S^{1}$ gives

\begin{equation}\label{2.10}
\| \nabla (P_{j} u \overline{P_{k} v}) \|_{L_{t,x}^{2}(\mathbf{R} \times \mathbf{R}^{d})}^{2} \lesssim 2^{(d - 1)k} 2^{j} \| P_{j} u \|_{L_{t}^{\infty} L_{x}^{2}(\mathbf{R} \times \mathbf{R}^{d})}^{2} \| P_{k} v \|_{L_{t}^{\infty} L_{x}^{2}(\mathbf{R} \times \mathbf{R}^{d})}^{2}.
\end{equation}

\noindent Thus, if $j \geq k + 10$, then by the Fourier support of $P_{j} u \cdot \overline{P_{k} v}$, and Bernstein's inequality,

\begin{equation}\label{2.11}
\| (P_{j} u \overline{P_{k} v}) \|_{L_{t,x}^{2}(\mathbf{R} \times \mathbf{R}^{d})}^{2} \lesssim 2^{(d - 1)k} 2^{-j} \| P_{j} u \|_{L_{t}^{\infty} L_{x}^{2}(\mathbf{R} \times \mathbf{R}^{d})}^{2} \| P_{k} v \|_{L_{t}^{\infty} L_{x}^{2}(\mathbf{R} \times \mathbf{R}^{d})}^{2}.
\end{equation}

\noindent Also, by theorem $\ref{t1.1}$, when $d \geq 3$, and $j - 10 \leq k \leq j$,

\begin{equation}\label{2.12}
\| (P_{j} u \overline{P_{k} v}) \|_{L_{t,x}^{2}(\mathbf{R} \times \mathbf{R}^{d})}^{2} \lesssim 2^{(d - 1)k} 2^{-j} \| P_{j} u \|_{L_{t}^{\infty} L_{x}^{2}(\mathbf{R} \times \mathbf{R}^{d})}^{2} \| P_{k} v \|_{L_{t}^{\infty} L_{x}^{2}(\mathbf{R} \times \mathbf{R}^{d})}^{2}.
\end{equation}

\noindent \textbf{Remark:} This will not work in dimensions $d = 1, 2$ because we can merely say that $P_{j} u \cdot \overline{P_{k} v}$ is supported on $|\xi| \lesssim 2^{k}$ when $j - 10 \leq k \leq j$, so when the power of $|\nabla|^{\frac{3 - d}{2}}$ is positive, then $(\ref{2.12})$ does not follow directly from theorem $\ref{t1.1}$ and the above analysis. It is for this technical reason that this paper only addresses dimensions $d \geq 3$.\vspace{5mm}

\noindent It is appropriate to point out that $(\ref{2.11})$ has been proven to be true in \cite{B2} using Fourier analysis. In fact, \cite{B2} and subsequent Fourier analytic proofs in higher dimensions all predate the bilinear virial arguments of \cite{PV}. However, in this case, it is quite useful to examine the bilinear estimates through the lens of the virial identity. In \cite{PV} the interaction Morawetz for theorem $\ref{t2.2}$ is

\begin{equation}\label{2.13}
\aligned
M(t) = \int |v(t,y)|^{2} \frac{(x - y)_{1}}{|(x - y)_{1}|} Im[\bar{u} \partial_{1} u](t,x) dx dy \\ + \int |u(t,y)|^{2} \frac{(x - y)_{1}}{|(x - y)_{1}|} Im[\bar{v} \partial_{1} v](t,x) dx dy.
\endaligned
\end{equation}

\noindent Averaging over all directions in $S^{d - 1}$, suppose $u$ is the solution to

\begin{equation}\label{2.14}
i u_{t} + \Delta u = \mathcal N_{1},
\end{equation}

\noindent and $v$ is the solution to

\begin{equation}\label{2.15}
i v_{t} + \Delta v = \mathcal N_{2}.
\end{equation}

\noindent For any $\mathbf{R}^{d}$, for some $c_{d}$, for any $y \in \mathbf{R}^{d}$,

\begin{equation}\label{2.15.1}
\int_{S^{d - 1}} \frac{x \cdot \omega}{|x \cdot \omega|} (y \cdot \omega) d\omega = c_{d} \frac{x \cdot y}{|x|}.
\end{equation}

\noindent Without loss of generality suppose $x = |x| e_{1}$. Then,

\begin{equation}\label{2.15.2}
\frac{1}{2} \frac{x \cdot \omega}{|x \cdot \omega|} y \cdot \omega + \frac{1}{2} \frac{x \cdot (-\omega)}{|x \cdot (-\omega)|} (y \cdot (-\omega)) = \frac{x_{1} \omega_{1}^{2}}{|x_{1}| |\omega_{1}|} y_{1},
\end{equation}

\noindent which implies $(\ref{2.15.1})$. Then by theorem $\ref{t2.2}$ and averaging over $\omega \in S^{d - 1}$, taking $\partial_{\omega} f = \omega \cdot \nabla f$, and $x_{\omega} = x \cdot \omega$,

\begin{equation}\label{2.16}
\aligned
\int_{S^{d - 1}} \int_{x_{\omega} = y_{\omega}} |\partial_{\omega}(u(t,y) \overline{v(t,x)})|^{2} dx dy dt d\omega \\
 \lesssim \| u \|_{L_{t}^{\infty} \dot{H}^{1/2}(\mathbf{R} \times \mathbf{R}^{d})}^{2} \| v \|_{L_{t}^{\infty} L_{x}^{2}(\mathbf{R} \times \mathbf{R}^{d})}^{2} +  \| v \|_{L_{t}^{\infty} \dot{H}^{1/2}(\mathbf{R} \times \mathbf{R}^{d})}^{2} \| u \|_{L_{t}^{\infty} L_{x}^{2}(\mathbf{R} \times \mathbf{R}^{d})}^{2} \\ + \int | v(t,y)|^{2} \frac{(x - y)}{|(x - y)|} \cdot Im[\overline{u} \nabla \mathcal N_{1}](t,x) dx dy dt \\
+ \int |v(t,y)|^{2} \frac{(x - y)}{|(x - y)|} \cdot Im[\overline{\mathcal N_{1}} \nabla u](t,x) dx dy dt \\ +  \int |u(t,y)|^{2} \frac{(x - y)}{|(x - y)|} \cdot Im[\overline{v} \nabla \mathcal N_{2}](t,x) dx dy dt \\ + \int |u(t,y)|^{2} \frac{(x - y)}{|(x - y)|} \cdot Im[\overline{\mathcal N_{2}} \nabla v](t,x) dx dy dt \\ + \int Re[\overline{v} \mathcal N_{2}](t,y) \frac{(x - y)}{|(x - y)|} \cdot Im[\bar{u} \nabla u](t,x) dx dy dt \\ + \int Re[\overline{u} \mathcal N_{1}](t,y) \frac{(x - y)}{|(x - y)|} \cdot Im[\bar{v} \nabla v](t,x) dx dy dt.
\endaligned
\end{equation}

\noindent Then by H{\"o}lder's inequality and averaging over $\omega \in S^{d - 1}$, for any $R > 0$,

\begin{equation}\label{2.16.1}
\aligned
\frac{1}{R} \int \int_{h \in \mathbf{R}^{d} : |h| \leq R} \int_{\mathbf{R}^{d}} |\nabla (u(t,x) \overline{v(t,x + h)})|^{2} dx dh dt \\
 \lesssim \| u \|_{L_{t}^{\infty} \dot{H}^{1/2}(\mathbf{R} \times \mathbf{R}^{d})}^{2} \| v \|_{L_{t}^{\infty} L_{x}^{2}(\mathbf{R} \times \mathbf{R}^{d})}^{2} +  \| v \|_{L_{t}^{\infty} \dot{H}^{1/2}(\mathbf{R} \times \mathbf{R}^{d})}^{2} \| u \|_{L_{t}^{\infty} L_{x}^{2}(\mathbf{R} \times \mathbf{R}^{d})}^{2} \\ + \int | v(t,y)|^{2} \frac{(x - y)}{|(x - y)|} \cdot Im[\overline{u} \nabla \mathcal N_{1}](t,x) dx dy dt \\
+ \int |v(t,y)|^{2} \frac{(x - y)}{|(x - y)|} \cdot Im[\overline{\mathcal N_{1}} \nabla u](t,x) dx dy dt \\ +  \int |u(t,y)|^{2} \frac{(x - y)}{|(x - y)|} \cdot Im[\overline{v} \nabla \mathcal N_{2}](t,x) dx dy dt \\ + \int |u(t,y)|^{2} \frac{(x - y)}{|(x - y)|} \cdot Im[\overline{\mathcal N_{2}} \nabla v](t,x) dx dy dt \\ + \int Re[\overline{v} \mathcal N_{2}](t,y) \frac{(x - y)}{|(x - y)|} \cdot Im[\bar{u} \nabla u](t,x) dx dy dt \\ + \int Re[\overline{u} \mathcal N_{1}](t,y) \frac{(x - y)}{|(x - y)|} \cdot Im[\bar{v} \nabla v](t,x) dx dy dt.
\endaligned
\end{equation}

\noindent Therefore, for $k \leq j$,

\begin{equation}\label{2.17}
\aligned
\| (P_{j} u \overline{P_{k} v}) \|_{L_{t,x}^{2}(\mathbf{R} \times \mathbf{R}^{d})}^{2} \lesssim 2^{(d - 1)k} 2^{-j} \| P_{j} u \|_{L_{t}^{\infty} L_{x}^{2}(\mathbf{R} \times \mathbf{R}^{d})}^{2} \| P_{k} v \|_{L_{t}^{\infty} L_{x}^{2}(\mathbf{R} \times \mathbf{R}^{d})}^{2} \\ + 2^{(d - 1)k} 2^{-2j} \int |P_{k} v(t,y)|^{2} \frac{(x - y)}{|(x - y)|} \cdot Re[\overline{P_{j} u} \nabla P_{j} \mathcal N_{1}](t,x) dx dy dt \\
+ 2^{(d - 1)k} 2^{-2j} \int |P_{k} v(t,y)|^{2} \frac{(x - y)}{|(x - y)|} \cdot Re[\overline{P_{j} \mathcal N_{1}} \nabla P_{j} u](t,x) dx dy dt \\ + 2^{(d - 1)k} 2^{-2j} \int |P_{j} u(t,y)|^{2} \frac{(x - y)}{|(x - y)|} \cdot Re[\overline{P_{k} v} \nabla P_{k} \mathcal N_{2}](t,x) dx dy dt \\ + 2^{(d - 1)k} 2^{-2j} \int |P_{j} u(t,y)|^{2} \frac{(x - y)}{|(x - y)|} \cdot Re[\overline{P_{k} \mathcal N_{2}} \nabla P_{k} v](t,x) dx dy dt \\ + 2^{(d - 1)k} 2^{-2j} \int Re[\overline{P_{k} v} P_{k} \mathcal N_{2}](t,y) \frac{(x - y)}{|(x - y)|} \cdot Re[P_{j} \bar{u} \nabla P_{j} u](t,x) dx dy dt \\ + 2^{(d - 1)k} 2^{-2j} \int Re[\overline{P_{j} u} P_{j} \mathcal N_{1}](t,y) \frac{(x - y)}{|(x - y)|} \cdot Re[P_{k} \bar{v} \nabla P_{k} v](t,x) dx dy dt.
\endaligned
\end{equation}

\noindent This is a useful estimate since it opens up the possibility of integrating by parts, and moving a derivative to a more advantageous position. This opportunity will prove to be quite useful, since we will follow the analysis in \cite{BIK}. That is, if $\phi$ be a solution to $(\ref{1.1})$, let

\begin{equation}\label{2.18}
\psi_{x} = \nabla \phi = (\psi_{1}, ..., \psi_{d}) = (\partial_{1} \phi, ..., \partial_{d} \phi).
\end{equation}

\noindent Then, choosing to work in the Coulomb gauge, for $m = 1, ..., d$,

\begin{equation}\label{2.19}
(i \partial_{t} + \Delta_{x}) \psi_{m} = -2 i \sum_{l = 1}^{d} A_{l} \cdot \partial_{l} \psi_{m} + (A_{0} + \sum_{l = 1}^{d} A_{l}^{2}) \psi_{m} - i \sum_{l = 1}^{d} Im(\psi_{m} \bar{\psi}_{l}) \psi_{l},
\end{equation}

\noindent where

\begin{equation}\label{2.20}
A_{l} = -\sum_{k = 1}^{d} \frac{\partial_{k}}{\Delta} Im(\psi_{l} \bar{\psi}_{k}),
\end{equation}

\noindent and

\begin{equation}\label{2.21}
A_{0} = \sum_{l, m = 1}^{d} \frac{\partial_{l} \partial_{m}}{\Delta} (Re(\bar{\psi}_{l} \psi_{m})) - \frac{1}{2} (\sum_{m = 1}^{d} \psi_{m} \bar{\psi}_{m}).
\end{equation}

\section{Growth of Besov norms}
Now we utilize frequency envelopes to control the growth of the Besov norm. The frequency envelopes were introduced in \cite{Taowm} for the study of wave maps. Frequency envelopes majorize the size of a Littlewood - Paley projection of a function, while at the same time smoothing out the differences in size at different frequency levels. Let $\delta > 0$ be a small constant, say $\delta = \frac{1}{4}$, and let

\begin{equation}\label{3.1}
\alpha_{j}(0) = \sup_{k} 2^{-\delta |j - k|} 2^{k(\frac{d - 2}{2})} \| P_{k} \psi_{x}(0) \|_{L^{2}(\mathbf{R}^{d})},
\end{equation}

\begin{equation}\label{3.2}
\alpha_{j} = \sup_{k} 2^{-\delta |j - k|} 2^{k(\frac{d - 2}{2})} \| P_{k} \psi_{x} \|_{L_{t}^{\infty} L_{x}^{2}(\mathbf{R} \times \mathbf{R}^{d})},
\end{equation}

\begin{equation}\label{3.3}
\beta_{j} = \sup_{k} 2^{-(\delta/2) |j - k|} 2^{k(\frac{d - 2}{4})} \| P_{k} \psi_{x} \|_{L_{t,x}^{4}(\mathbf{R} \times \mathbf{R}^{d})},
\end{equation}

\noindent and

\begin{equation}\label{3.4}
\aligned
\gamma_{j} = \sup_{k} 2^{-\delta |j - k|} 2^{k(\frac{d - 2}{2})} (\sup_{l \leq k} 2^{\frac{1}{2}(k - l)} \| (P_{k} \psi_{x})(P_{l} \psi_{x}) \|_{L_{t,x}^{2}(\mathbf{R} \times \mathbf{R}^{d})}) \\ + \sup_{k} 2^{-\delta |j - k|} 2^{k(\frac{d - 2}{2})} (\sup_{l \leq k} \sup_{R > 0} \frac{1}{R} 2^{k} 2^{(d - 2)l} \int_{h \in \mathbf{R}^{d}, |h| \leq R} |P_{k} \psi_{x}(t,x)|^{2} |P_{l} \psi_{x}(t, x + h)|^{2} dx dy dt)^{1/2}.
\endaligned
\end{equation}

\noindent \textbf{Remark:} From the definition,

\begin{equation}\label{3.5}
\beta_{j}^{2} \lesssim \gamma_{j}.
\end{equation}

\noindent Also observe that for any $j$,

\begin{equation}
\| P_{j} \psi_{x}(0) \|_{\dot{H}^{\frac{d - 2}{2}}} \leq \alpha_{j}(0),
\end{equation}

\noindent and by Young's inequality,

\begin{equation}
\| \psi_{x}(0) \|_{\dot{H}^{\frac{d - 2}{2}}}^{2} \sim \sum_{j} \alpha_{j}(0)^{2}.
\end{equation}

\noindent Now if $\langle f , g \rangle$ is the inner product

\begin{equation}\label{3.6}
Re \int f(x) \overline{g(x)} dx,
\end{equation}

\noindent then

\begin{equation}\label{3.7}
\| P_{j} \psi_{x} \|_{L_{x}^{2}(\mathbf{R}^{d})}^{2} = \langle P_{j} \psi_{x}, P_{j} \psi_{x} \rangle,
\end{equation}

\noindent so by $(\ref{2.19})$,

\begin{equation}\label{3.8}
\aligned
\frac{d}{dt} \langle P_{j} \psi_{x}, P_{j} \psi_{x} \rangle = -4 \langle P_{j}(\sum_{l = 1}^{d} A_{l} \partial_{l} \psi_{x}), P_{j} \psi_{x} \rangle \\ - 2 \langle i P_{j}((A_{0} + \sum_{l = 1}^{d} A_{l}^{2}) \psi_{x}), P_{j} \psi_{x} \rangle - 2 \langle P_{j}(\sum_{l = 1}^{d} Im(\psi_{x} \bar{\psi}_{l}) \psi_{l}), P_{j} \psi_{x} \rangle.
\endaligned
\end{equation}

\noindent Now by Fourier support arguments,

\begin{equation}\label{3.9}
\aligned
\sup_{k} 2^{-2 \delta |j - k|} 2^{k(d - 2)} \int_{\mathbf{R}} |\langle P_{k}(\sum_{l = 1}^{d} Im(\psi_{x} \bar{\psi}_{l}) \psi_{l}), P_{k} \psi_{x} \rangle| dt \\ \lesssim \sup_{k} 2^{-2 \delta |j - k|} 2^{k(d - 2)} \| P_{k} \psi_{x} \|_{L_{t,x}^{4}} \| P_{\geq k - 10} \psi_{x} \|_{L_{t,x}^{4}}^{3} \\ + \sup_{k} 2^{-2\delta |j - k|} 2^{k(d - 2)} \| (P_{k} \psi_{x})(P_{\leq k - 10} \psi_{x}) \|_{L_{t,x}^{2}} \| (P_{k - 10 \leq \cdot \leq k + 10} \psi_{x})(P_{\leq k - 10} \psi_{x}) \|_{L_{t,x}^{2}} \\ \lesssim \beta_{j}^{4} + \gamma_{j}^{2} \lesssim \gamma_{j}^{2}.
\endaligned
\end{equation}

\noindent By a similar calculation,

\begin{equation}\label{3.10}
\sup_{k} 2^{-2 \delta |j - k|} 2^{k(d - 2)} \int_{\mathbf{R}} |\langle i P_{k}(-\frac{1}{2} (\sum_{m = 1}^{d} \psi_{m} \bar{\psi}_{m}) \psi_{x}), P_{k} \psi_{x} \rangle| dt \lesssim \beta_{j}^{4} + \gamma_{j}^{2} \lesssim \gamma_{j}^{2}.
\end{equation}

\noindent Now split

\begin{equation}\label{3.11}
\aligned
\frac{\partial_{l} \partial_{m}}{\Delta} Re(\bar{\psi}_{l} \psi_{m}) =  \frac{\partial_{l} \partial_{m}}{\Delta} Re((P_{\geq k - 20} \bar{\psi}_{l})(P_{\geq k - 20} \psi_{m})) \\ + 2 \frac{\partial_{l} \partial_{m}}{\Delta} Re((P_{\geq k - 20} \bar{\psi}_{l})(P_{\leq k - 20} \psi_{m})) + \frac{\partial_{l} \partial_{m}}{\Delta} Re((P_{\leq k - 20} \bar{\psi}_{l})(P_{\leq k - 20} \psi_{m})).
\endaligned
\end{equation}

\noindent Because $\frac{\partial_{l} \partial_{m}}{\Delta}$ is a bounded Fourier multiplier,

\begin{equation}\label{3.12}
\aligned
\sup_{k} 2^{-2 \delta |j - k|} 2^{k(d - 2)} \int |\langle \frac{\partial_{l} \partial_{m}}{\Delta} (Re(P_{\geq k - 20} \bar{\psi}_{l}) (P_{\geq k - 20} \psi_{m})) \psi_{x}, P_{k} \psi_{x} \rangle| dt \\
\lesssim \sup_{k} 2^{-2 \delta |j - k|} 2^{k(d - 2)} \| P_{\geq k - 20} \psi_{x} \|_{L_{t,x}^{4}}^{2} \| P_{k} \psi_{x} \|_{L_{t,x}^{4}} \| P_{\geq k - 10} \psi_{x} \|_{L_{t,x}^{4}} \\ + \sup_{k} 2^{-2 \delta |j - k|} 2^{k(d - 2)} \| P_{\geq k - 20} \psi_{x} \|_{L_{t,x}^{4}}^{2} \| (P_{k} \psi_{x})(P_{\leq k - 10} \psi_{x}) \|_{L_{t,x}^{2}} \lesssim \beta_{j}^{4} + \gamma_{j}^{2}.
\endaligned
\end{equation}

\noindent Meanwhile, because $(P_{\leq k - 20} \bar{\psi}_{l})(P_{\leq k - 20} \psi_{m})$ is supported on $|\xi| \leq 2^{k - 15}$, by $(\ref{2.7})$ we compute

\begin{equation}\label{3.13}
\aligned
\sup_{k} 2^{-2 \delta |j - k|} 2^{(d - 2)k} \int \sum_{l \leq k - 15} \int_{\mathbf{R}^{d}} {2^{ld}} |\check{\psi}(2^{l}(x - y))| |P_{l \leq \cdot \leq k - 20} \psi_{x}(t,y)|^{2} \\ \times |P_{k} \psi_{x}(t,x)| |P_{k - 10 \leq \cdot \leq k + 10} \psi_{x}(t,x)| dx dy dt \\
\lesssim \gamma_{j}^{2} \sum_{l \leq k - 20} 2^{(d - 1)l} \sum_{l \leq m \leq k - 20} 2^{-(d - 2) m} 2^{-(d - 2)k} 2^{-k} \lesssim \gamma_{j}^{2}.
\endaligned
\end{equation}

\noindent Also by $(\ref{2.7})$,

\begin{equation}\label{3.14}
\aligned
\sup_{k} 2^{-2 \delta |j - k|} 2^{(d - 2)k} \int \sum_{l \leq k - 15} \int_{\mathbf{R}^{d}} {2^{ld}} |\check{\psi}(2^{l}(x - y))| |P_{\leq l} \psi_{x}(t,y)|^{2} \\ \times |P_{k} \psi_{x}(t,x)| |P_{k - 10 \leq \cdot \leq k + 10} \psi_{x}(t,x)| dx dy dt \\
\lesssim \sup_{k} 2^{-2 \delta |j - k|} 2^{(d - 2)k} \sum_{l \leq k - 15} \| (P_{k} \psi_{x})(P_{\leq l} \psi_{x}) \|_{L_{t,x}^{2}} \| (P_{k - 10 \leq \cdot \leq k + 10} \psi_{x})(P_{\leq l} \psi_{x}) \|_{L_{t,x}^{2}} \lesssim \gamma_{j}^{2}.
\endaligned
\end{equation}

\noindent Therefore, combining $(\ref{3.10})$ - $(\ref{3.14})$ with an interpolation argument proves

\begin{equation}\label{3.15}
\sup_{k} 2^{-2 \delta |j - k|} 2^{k(d - 2)} |\langle i P_{k}((A_{0}) \psi_{x}), P_{k} \psi_{x} \rangle| dt \lesssim \gamma_{j}^{2}.
\end{equation}

\noindent Next, integrating by parts, since $\sum_{l = 1}^{d} \partial_{l} A_{l} = 0$

\begin{equation}\label{3.16}
-4 \langle \sum_{l = 1}^{d} A_{l} \cdot \partial_{l} (P_{k} \psi_{x}), P_{k} \psi_{x} \rangle = 0.
\end{equation}

\noindent To estimate

\begin{equation}\label{3.16.1}
-4 \langle \sum_{l = 1}^{d} [P_{k}, A_{l}] \partial_{l} \psi_{x}, P_{k} \psi_{x} \rangle,
\end{equation}

\noindent observe that if $m$ is a Fourier multiplier satisfying $|\nabla m(\xi)| \lesssim \frac{1}{|\xi|}$ and $|m(\xi)| \lesssim 1$, then by the fundamental theorem of calculus, if $|\eta| << |\xi|$,

\begin{equation}\label{3.17}
|\xi| |m(\xi + \eta) \hat{f}(\eta) \hat{g}(\xi) - m(\xi) \hat{f}(\eta) \hat{g}(\xi)| \lesssim |\eta| |\hat{f}(\eta)| |\hat{g}(\xi)|.
\end{equation}

\noindent If $|\eta| \gtrsim |\xi|$,

\begin{equation}\label{3.18}
|\xi| |m(\xi + \eta) \hat{f}(\eta) \hat{g}(\xi)| + |\xi| |m(\xi) \hat{f}(\eta) \hat{g}(\xi)| \lesssim |\eta| |\hat{f}(\eta)| |\hat{g}(\xi)|.
\end{equation}

\noindent Therefore, by $(\ref{2.20})$

\begin{equation}\label{3.19}
2^{-2 \delta |j - k|} 2^{k(d - 2)} \int |\langle [P_{k}, A_{l}] \partial_{l} \psi_{x}, P_{k} \psi_{x} \rangle| dt
\end{equation}

\noindent may be estimated in a manner identical to the argument leading up to $(\ref{3.15})$.\vspace{5mm}

\noindent It finally remains to estimate

\begin{equation}\label{3.20}
\sup_{k} 2^{-2 \delta |j - k|} 2^{(d - 2)k} \int |\langle i P_{k} (\sum_{l = 1}^{d} A_{l}^{2} \psi_{x}), P_{k} \psi_{x} \rangle| dt.
\end{equation}

\noindent By the Sobolev embedding theorem and Bernstein's inequality,

\begin{equation}\label{3.21}
\aligned
\| \psi_{x} \|_{L_{x}^{d}(\mathbf{R}^{d})}^{d} \lesssim \sum_{j_{1} \leq ... \leq j_{d}} \| P_{j_{1}} \|_{L_{x}^{\infty}} \cdots \| P_{j_{d - 2}} \|_{L_{x}^{\infty}} \| P_{j_{d - 1}} \psi_{x} \|_{L_{x}^{2}} \| P_{j_{d}} \psi_{x} \|_{L_{x}^{2}} \\
\lesssim \sum_{j_{1} \leq ... \leq j_{d}} \| P_{j_{1}} \psi_{x} \|_{\dot{H}_{x}^{\frac{d - 2}{2}}} \cdots \| P_{j_{d}} \|_{\dot{H}_{x}^{\frac{d - 2}{2}}} 2^{j_{1}} \cdots 2^{j_{d - 2}} 2^{-j_{d - 1} (\frac{d - 2}{2})} 2^{-j_{d} (\frac{d - 2}{2})} \\
\lesssim \| \psi_{x} \|_{\dot{H}_{x}^{\frac{d - 2}{2}}}^{2} (\sup_{j} \| P_{j} \psi_{x} \|_{\dot{H}^{\frac{d - 2}{2}}})^{d - 2}.
\endaligned
\end{equation}

\noindent Next, by the Sobolev embedding theorem and $(\ref{3.21})$,

\begin{equation}\label{3.22}
\| P_{k}(-\sum_{k = 1}^{d} \frac{\partial_{k}}{\Delta} Im(\psi_{l} \bar{\psi}_{k})) \|_{L_{x}^{\infty}} \lesssim 2^{k} \| \psi_{x} \|_{L_{x}^{d}}^{2} \lesssim 2^{k} \| \psi_{x} \|_{\dot{H}_{x}^{\frac{d - 2}{2}}}^{\frac{4}{d}} (\sup_{j} \| P_{j} \psi_{x} \|_{\dot{H}^{\frac{d - 2}{2}}})^{\frac{2(d - 2)}{d}}.
\end{equation}

\noindent Therefore, by Fourier support properties,

\begin{equation}\label{3.23}
\aligned
\sum_{n \leq m \leq k - 10} \| |P_{k}(P_{m} A_{x} \cdot P_{n} A_{x} \cdot \psi_{x})| |P_{k} \psi_{x}| \|_{L_{t,x}^{1}} \\
\lesssim \| \psi_{x} \|_{\dot{H}_{x}^{\frac{d - 2}{2}}}^{\frac{4}{d}} (\sup_{j} \| P_{j} \psi_{x} \|_{\dot{H}^{\frac{d - 2}{2}}})^{\frac{2(d - 2)}{d}} \sum_{m} \| |P_{k} (2^{m} P_{m} A_{x} \cdot \psi_{x})| |P_{k} \psi_{x}| \|_{L_{t,x}^{1}}.
\endaligned
\end{equation}

\noindent Now since the Fourier multiplier $2^{m} \frac{\partial_{k}}{\Delta} P_{m}$ is uniformly bounded, by an argument similar to $(\ref{3.11})$ - $(\ref{3.15})$,

\begin{equation}\label{3.24}
\sup_{k} 2^{-2 \delta |j - k|} 2^{(d - 2)k} (\ref{3.23}) \lesssim \| \psi_{x} \|_{\dot{H}_{x}^{\frac{d - 2}{2}}}^{\frac{4}{d}} (\sup_{j} \| P_{j} \psi_{x} \|_{\dot{H}^{\frac{d - 2}{2}}})^{\frac{2(d - 2)}{d}}(\beta_{j}^{4} + \gamma_{j}^{2}) \lesssim \gamma_{j}^{2}.
\end{equation}

\noindent On the other hand, if $m \geq k - 10$, then by Fourier support arguments and the Sobolev embedding theorem,

\begin{equation}\label{3.25}
\| P_{n} A_{x} \|_{L_{t,x}^{4}} \lesssim \| P_{\geq k - 15} \psi_{x} \|_{L_{t,x}^{4}} \| \psi_{x} \|_{L_{x}^{d}},
\end{equation}

\noindent and therefore,

\begin{equation}\label{3.26}
\aligned
\sup_{k} 2^{-2 \delta |j - k|} 2^{(d - 2)k} \int |\langle i P_{k} (\sum_{l = 1}^{d} A_{l}^{2} \psi_{x}), P_{k} \psi_{x} \rangle| dt \\ \lesssim \| \psi_{x} \|_{\dot{H}_{x}^{\frac{d - 2}{2}}}^{\frac{4}{d}} (\sup_{j} \| P_{j} \psi_{x} \|_{\dot{H}^{\frac{d - 2}{2}}})^{\frac{2(d - 2)}{d}}(\beta_{j}^{4} + \gamma_{j}^{2}).
\endaligned
\end{equation}

\noindent In conclusion, we have proved

\begin{theorem}\label{t3.1}
The frequency envelope $\alpha$ has the bounds

\begin{equation}\label{3.27}
\aligned
\alpha_{j} \lesssim \alpha_{j}(0) + \gamma_{j} (1 + \| \psi_{x} \|_{\dot{H}_{x}^{\frac{d - 2}{2}}}^{\frac{2}{d}} (\sup_{j} \| P_{j} \psi_{x} \|_{\dot{H}^{\frac{d - 2}{2}}})^{\frac{d - 2}{d}}) \\
\lesssim \alpha_{j}(0) + \gamma_{j}(1 + (\sum_{j} \alpha_{j}^{2})^{2/d} (\sup_{j} \alpha_{j}^{2})^{\frac{d - 2}{d}}).
\endaligned
\end{equation}
\end{theorem}

\section{Bilinear estimates}
Now we use the interaction Morawetz estimates to prove some bilinear estimates.

\begin{proposition}\label{p4.1}
For any $j$,

\begin{equation}\label{4.1}
\aligned
\gamma_{j}^{2} \lesssim (\sup_{k} \alpha_{k})^{2} \alpha_{j}^{2} + \gamma_{j}^{2} (\sup_{k} \alpha_{k})^{2} (1 + (\sum_{k} \alpha_{k}^{2})^{2/d} (\sup_{k} \alpha_{k})^{\frac{2(d - 2)}{d}}) \\
+ \alpha_{j}^{2} (\sup_{k} \gamma_{k})^{2} (1 + (\sum_{k} \alpha_{k}^{2})^{2/d} (\sup_{k} \alpha_{k})^{\frac{2(d - 2)}{d}}) .
\endaligned
\end{equation}
\end{proposition}

\noindent \emph{Proof:} Take $(\ref{2.17})$ and set $\mathcal N_{1} = \mathcal N_{2} = \mathcal N$, $u = v = \psi_{x}$, where

\begin{equation}\label{4.2}
\mathcal N = - i \sum_{l = 1}^{d} Im(\psi_{x} \bar{\psi}_{l}) \psi_{l} + A_{0} \psi_{x} + \sum_{l = 1}^{d} A_{l}^{2} \psi_{x}  - 2 i \sum_{l = 1}^{d} A_{l} \cdot \partial_{l} \psi_{x} = \mathcal N^{(1)} + \mathcal N^{(2)} + \mathcal N^{(3)} + \mathcal N^{(4)}.
\end{equation}

\noindent We will call the first line on the right hand side of $(\ref{2.17})$ the main term, and the other lines the error. Then we label the contribution to the error from $\mathcal N^{(1)}$, $\mathcal E^{(1)}$, and so on. The main term may be easily estimated.

\begin{equation}\label{4.3}
2^{(d - 2)l} 2^{-l} \| P_{j} \psi_{x} \|_{L_{t}^{\infty} L_{x}^{2}}^{2} \sup_{k \leq l} (2^{l - k} 2^{(d - 1)k} \| P_{k} \psi_{x} \|_{L_{t}^{\infty} L_{x}^{2}}^{2}) \lesssim \alpha_{j}^{2} (\sup_{k} \alpha_{k})^{2}.
\end{equation}

\noindent Computing $\mathcal E^{(1)}$, $(\ref{3.9})$ implies

\begin{equation}\label{4.4}
\aligned
2^{l - k} 2^{(d - 2)l} \mathcal E^{(1)} \lesssim 2^{(d - 2)l} 2^{(d - 2)k} \| |P_{l} \mathcal N^{(1)}| |P_{l} \psi_{x}| \|_{L_{t,x}^{1}} \| P_{k} \psi_{x} \|_{L_{t}^{\infty} L_{x}^{2}}^{2} \\ + 2^{(d - 2)k} 2^{(d - 2)l} 2^{k - l} \| |P_{k} \mathcal N^{(1)}| |P_{k} \psi_{x}| \|_{L_{t,x}^{1}} \| P_{l} \psi_{x} \|_{L_{t}^{\infty} L_{x}^{2}}^{2} \\
\lesssim 2^{2\delta |j - l|} [\gamma_{j}^{2} (\sup_{k} \alpha_{k})^{2} + \alpha_{j}^{2} (\sup_{k} \gamma_{k})^{2}].
\endaligned
\end{equation}

\noindent Similarly, by $(\ref{3.10})$ - $(\ref{3.16})$,

\begin{equation}\label{4.5}
2^{l - k} 2^{(d - 2)l} \mathcal E^{(2)} \lesssim 2^{2\delta |j - l|} [\gamma_{j}^{2} (\sup_{k} \alpha_{k})^{2} + \alpha_{j}^{2} (\sup_{k} \gamma_{k})^{2}],
\end{equation}

\noindent and by $(\ref{3.20})$ - $(\ref{3.26})$,

\begin{equation}\label{4.6}
2^{l - k} 2^{(d - 2)l} \mathcal E^{(3)} \lesssim 2^{2\delta |j - l|} [\gamma_{j}^{2} (\sup_{k} \alpha_{k})^{2} + \alpha_{j}^{2} (\sup_{k} \gamma_{k})^{2}] (1 + \| \psi_{x} \|_{L_{t}^{\infty} L_{x}^{d}}^{2}).
\end{equation}

\noindent To estimate the integral of $\mathcal E^{(4)}$, split

\begin{equation}\label{4.7}
-4 P_{k} (A_{x} \cdot \nabla_{x} \psi_{x}) = -4 A_{x} \cdot \nabla_{x} (P_{k} \psi_{x}) - 4 [P_{k}, A_{x}] \cdot \nabla_{x} \psi_{x} = \mathcal N^{(4)}_{1} + \mathcal N^{(4)}_{2},
\end{equation}

\noindent and split $\mathcal E^{(4)} = \mathcal E^{(4)}_{1} + \mathcal E^{(4)}_{2}$. Then by $(\ref{3.17})$ - $(\ref{3.19})$,

\begin{equation}\label{4.8}
2^{l - k} 2^{(d - 2)l} \mathcal E^{(4)}_{2} \lesssim 2^{2\delta |j - l|} [\gamma_{j}^{2} (\sup_{k} \alpha_{k})^{2} + \alpha_{j}^{2} (\sup_{k} \gamma_{k})^{2}].
\end{equation}

\noindent To estimate the contribution of $\mathcal E^{(4)}_{1}$ we again integrate by parts. First,

\begin{equation}\label{4.9}
-4 Re[\overline{P_{k} \psi_{x}} (A_{x} \cdot \nabla_{x} (P_{k} \psi_{x}))] = -2 Re[A_{x} \cdot \nabla_{x} (|P_{k} \psi_{x}|^{2})].
\end{equation}

\noindent Then integrating by parts,

\begin{equation}\label{4.10}
2^{-2l} 2^{(d - 1)k} \int Re[\overline{P_{k} \psi_{x}} (A_{x} \cdot \nabla_{x} (P_{k} \psi_{x}))](t,y) \frac{(x - y)}{|(x - y)|} \cdot Re[P_{l} \bar{\psi}_{x} \nabla P_{l} \psi_{x}](t,x) dx dy dt
\end{equation}

\begin{equation}\label{4.11}
\lesssim 2^{-2l} 2^{(d - 1)k}\int |P_{k} \psi_{x}(t,y)|^{2} |A_{x}(t,y)| \frac{1}{|x - y|} |\nabla P_{l} \psi_{x}(t,x)| |P_{l} \psi_{x}(t,x)| dx dy dt.
\end{equation}

\noindent Now recall the definition of $A_{x}$ in $(\ref{2.20})$. Then we can split

\begin{equation}\label{4.12}
\aligned
A_{x} = \sum_{m = 1}^{d} \frac{\partial_{m}}{\Delta} Im(\psi_{x} \bar{\psi}_{m}) = \sum_{m = 1}^{d} \frac{\partial_{m}}{\Delta} Im((P_{\leq k - 10} \psi_{x}) (P_{\leq k - 10} \bar{\psi}_{m})) \\
+ 2 \sum_{m = 1}^{d} \frac{\partial_{m}}{\Delta} Im((P_{\leq k - 10} \psi_{x}) (P_{\geq k - 10} \bar{\psi}_{m})) + \sum_{m = 1}^{d} \frac{\partial_{m}}{\Delta} Im((P_{\geq k - 10} \psi_{x}) (P_{\geq k - 10} \bar{\psi}_{m})) \\ = A_{x}^{(1)} + A_{x}^{(2)} + A_{x}^{(3)}.
\endaligned
\end{equation}

\noindent Now by definition of $\gamma_{k}$ and $\beta_{k}$, if $\delta < \frac{1}{4}$, then by the Sobolev embedding theorem,

\begin{equation}\label{4.13}
\| A_{x}^{(2)} \|_{L_{t}^{2} L_{x}^{\frac{2d}{d - 2}}} + \| A_{x}^{(3)} \|_{L_{t}^{2} L_{x}^{\frac{2d}{d - 2}}} \lesssim \| |P_{\leq k - 10} \psi_{x}| |P_{\geq k - 10} \psi_{x} \|_{L_{t,x}^{2}} + \| P_{\geq k - 10} \psi_{x} \|_{L_{t,x}^{4}} \lesssim \gamma_{k} + \beta_{k}^{2}.
\end{equation}

\noindent Therefore,

\begin{equation}\label{4.14}
\aligned
2^{-2l} \sup_{k} (2^{l - k} 2^{(d - 1)k}\int |P_{k} \psi_{x}(t,y)|^{2} |A_{x}^{(2)} (t,y)| \frac{1}{|x - y|} |\nabla P_{l} \psi_{x}(t,x)| |P_{l} \psi_{x}(t,x)| dx dy dt) \\
+ 2^{-2l} \sup_{k} (2^{l - k} 2^{(d - 1)k}\int |P_{k} \psi_{x}(t,y)|^{2} |A_{x}^{(3)} (t,y)| \frac{1}{|x - y|} |\nabla P_{l} \psi_{x}(t,x)| |P_{l} \psi_{x}(t,x)| dx dy dt) \\
\lesssim 2^{2 \delta |j - l|} 2^{-j(d - 2)} (\sup_{k} \beta_{k})^{2} (\sup_{k} \gamma_{k} + (\sup_{k} \beta_{k})^{2}) \alpha_{l}^{2} \lesssim 2^{2 \delta |j - l|} 2^{-l(d - 2)} ((\sup_{k} \gamma_{k})^{2} + (\sup_{k} \beta_{k})^{4}) \alpha_{j}^{2}.
\endaligned
\end{equation}

\noindent Now by Fourier support properties,

\begin{equation}\label{4.15}
\sum_{m = 1}^{d} \frac{\partial_{m}}{\Delta} Im((P_{\leq k - 10} \psi_{x}) (P_{\leq k - 10} \bar{\psi}_{m}))
\end{equation}

\noindent is supported on $|\xi| \leq 2^{k - 5}$, so

\begin{equation}\label{4.16}
A_{x}^{(1)} = \sum_{m \leq k - 5} P_{m} A_{x}^{(1)}.
\end{equation}

\noindent The Fourier multiplier $P_{m} (\frac{\partial_{k}}{\Delta})$ has size $\sim 2^{-m}$, so by $(\ref{3.13})$ and $(\ref{3.14})$,

\begin{equation}\label{4.17}
\aligned
2^{-2l} \sup_{k} (2^{l - k} 2^{(d - 1)k} \sum_{m \leq k - 5} \int_{|x - y| \geq 2^{-k} 2^{\frac{1}{d}(k - m)}} |P_{k} \psi_{x}(t,y)|^{2} |P_{m} A_{x}^{(1)} (t,y)| \\ \times \frac{1}{|x - y|} |\nabla P_{l} \psi_{x}(t,x)| |P_{l} \psi_{x}(t,x)| dx dy dt) \lesssim 2^{-l(d - 2)} 2^{2 \delta |j - l|} \alpha_{j}^{2} (\sup_{k} \gamma_{k})^{2}.
\endaligned
\end{equation}

\noindent Also, by the Sobolev embedding theorem and Bernstein's inequality,

\begin{equation}\label{4.18}
\| P_{m} A_{x}^{(1)} \|_{L_{x}^{\infty}} \lesssim 2^{-m} \| P_{\leq m} \psi_{x} \|_{L_{x}^{\infty}}^{2} + 2^{(d - 1)m} \| P_{\geq m} \psi_{x} \|_{L_{x}^{2}}^{2} \lesssim 2^{m} (\sup_{k} \alpha_{k})^{2}.
\end{equation}

\noindent Therefore, by H{\"o}lder's inequality,

\begin{equation}\label{4.19}
\aligned
2^{-2l} \sup_{k} (2^{l - k} 2^{(d - 1)k} \sum_{m \leq k - 5} \int_{|x - y| \leq 2^{-k} 2^{\frac{1}{d}(k - m)}} |P_{k} \psi_{x}(t,y)|^{2} |P_{m} A_{x}^{(1)} (t,y)| \\ \times \frac{1}{|x - y|} |\nabla P_{l} \psi_{x}(t,x)| |P_{l} \psi_{x}(t,x)| dx dy dt) \\ \lesssim  \sup_{k} \| (P_{l} \psi_{x})(P_{k} \psi_{x}) \|_{L_{t,x}^{2}}^{2} (\sup_{k} \alpha_{k}^{2}) \lesssim 2^{2 \delta |j - l|} 2^{-l(d - 2)} \gamma_{j}^{2} (\sup_{k} \alpha_{k})^{2}.
\endaligned
\end{equation}

\noindent Also,

\begin{equation}\label{4.20}
-4 Re[\overline{P_{l} \psi_{x}} (A_{x} \cdot \nabla_{x} (P_{l} \psi_{x}))] = -2 Re[A_{x} \cdot \nabla_{x} (|P_{l} \psi_{x}|^{2})],
\end{equation}

\noindent so following the analysis in $(\ref{4.9})$ - $(\ref{4.19})$ with $l$ and $k$ swapped proves

\begin{equation}\label{4.21}
\aligned
2^{-2l} \sup_{k} 2^{(d - 1)k} \int Re[\overline{P_{l} \psi_{x}} (A_{x} \cdot \nabla_{x} (P_{l} \psi_{x}))](t,y) \frac{(x - y)}{|(x - y)|} \cdot Re[P_{k} \bar{\psi}_{x} \nabla P_{k} \psi_{x}](t,x) dx dy dt \\
\lesssim 2^{-l(d - 2)} 2^{2 \delta |j - l|} \gamma_{j}^{2} (\sup_{k} \alpha_{k})^{2}.
\endaligned
\end{equation}

\noindent This finally proves proposition $\ref{p4.1}$. $\Box$\vspace{5mm}

\noindent \emph{Proof of theorem $\ref{t1.3}$:} Theorem $\ref{t1.3}$ may be proved by a bootstrap argument. We have a local existence result.

\begin{proposition}\label{p4.2}
Assume $\phi_{0} \in H_{Q}^{\infty}$. Then there is $T_{\sigma_{0}} = T(\| \phi_{0} \|_{H_{Q}^{\sigma_{0}}}) > 0$ and a solution $\phi \in C( [-T_{\sigma_{0}}, T_{\sigma_{0}}; H_{Q}^{\infty})$ of the initial value problem $(\ref{1.1})$. Additionally, $T_{\sigma_{0}}$ can be chosen so that

\begin{equation}\label{4.22}
\sup_{t \in [-T_{\sigma_{0}}, T_{\sigma_{0}}]} \| \phi(t) \|_{H_{Q}^{\sigma_{0}}} \leq C(\| \phi_{0} \|_{H_{Q}^{\sigma_{0}}}),
\end{equation}

\noindent and for $\sigma \in [\sigma_{0}, \infty) \cap \mathbf{Z}$,

\begin{equation}\label{4.23}
\sup_{t \in [-T_{\sigma_{0}}, T_{\sigma_{0}}]} \| \phi(t) \|_{H_{Q}^{\sigma}} \leq C(\sigma, \| \phi_{0} \|_{H_{Q}^{\sigma}}).
\end{equation}
\end{proposition}

\noindent \emph{Proof:} See \cite{KLPST}. $\Box$\vspace{5mm}

\noindent Then for $\phi_{0} \in H_{Q}^{\sigma_{0}}$, for $\sigma_{0}$ sufficiently large, and $T(\| \phi_{0} \|_{H_{Q}^{\sigma_{0}}}) > 0$ sufficiently small, we have that

\begin{equation}\label{4.24}
\sum_{k} \gamma_{k} \lesssim \epsilon,
\end{equation}

\noindent for $\epsilon(\| \phi_{0} \|_{\dot{H}^{d/2}})$ sufficiently small. Then by theorem $\ref{t3.1}$ and proposition $\ref{p4.1}$,

\begin{equation}\label{4.25}
\sum_{k} \alpha_{k}^{2} \lesssim \| \phi_{0} \|_{\dot{H}^{d/2}}^{2},
\end{equation}

\noindent and

\begin{equation}\label{4.26}
\sup_{k} \alpha_{k} + \sup_{k} \gamma_{k} \lesssim \epsilon.
\end{equation}

\noindent Again by theorem $\ref{t3.1}$ and proposition $\ref{p4.1}$,

\begin{equation}\label{4.27}
\alpha_{j} \lesssim \alpha_{j}(0) + \gamma_{j},
\end{equation}

\noindent and

\begin{equation}\label{4.28}
\gamma_{j} \lesssim \epsilon \alpha_{j} + \epsilon \gamma_{j}.
\end{equation}

\noindent Combining $(\ref{4.27})$ and $(\ref{4.28})$,

\begin{equation}\label{4.29}
\alpha_{j} \lesssim \alpha_{j}(0), \hspace{5mm} \text{and} \hspace{5mm} \gamma_{j} \lesssim \epsilon \alpha_{j}.
\end{equation}

\noindent But we could then replace $\alpha_{j}$ and $\gamma_{j}$ by $\tilde{\alpha}_{j}$ and $\tilde{\gamma}_{j}$, with

\begin{equation}\label{4.30}
\tilde{\alpha}_{j} = \sup_{k} 2^{-\delta |j - k|} 2^{k(\sigma_{0} - 1)} \| P_{k} \psi_{x} \|_{L_{t}^{\infty} L_{x}^{2}},
\end{equation}

\noindent and

\begin{equation}\label{4.31}
\tilde{\gamma}_{j} = \sup_{k} 2^{-\delta |j - k|} 2^{k(\sigma_{0} - 1)} \sup_{l} (2^{k - l} \| (P_{k} \psi_{x})(P_{l} \psi_{x}) \|_{L_{t,x}^{2}}
\end{equation}

\noindent Following the arguments proving theorem $\ref{t3.1}$ and proposition $\ref{p4.1}$, it is possible to prove that

\begin{equation}\label{4.32}
\tilde{\alpha}_{j} \lesssim \tilde{\alpha}_{j}(0) + \epsilon \tilde{\gamma}_{j},
\end{equation}

\noindent and

\begin{equation}\label{4.33}
\tilde{\gamma}_{j} \lesssim \epsilon \tilde{\alpha}_{j} + \epsilon \tilde{\gamma}_{j},
\end{equation}

\noindent and therefore

\begin{equation}\label{4.34}
\tilde{\alpha}_{j} \lesssim \tilde{\alpha}_{j}(0).
\end{equation}

\noindent In particular, this implies that for any $t \in [-T_{\sigma_{0}}, T_{\sigma_{0}}]$,

\begin{equation}\label{4.35}
\| \phi(t) \|_{\dot{H}^{\sigma_{0}}} \lesssim \| \phi_{0} \|_{\dot{H}^{\sigma}}.
\end{equation}

\noindent Then we may extend the interval of existence a little farther, but still maintaining $(\ref{4.24})$, which then implies that $(\ref{4.35})$ holds. We can iterate this argument, and the implicit constants in $(\ref{4.29})$ and $(\ref{4.35})$ are uniformly bounded. Making a standard bootstrap argument then implies theorem $\ref{t1.3}$. $\Box$

\nocite*
\bibliographystyle{plain}

\begin{thebibliography}{[00]}

\bibitem{B}
	\newblock I. Bejenaru,
	\newblock ``Global results for Schr{\"o}dinger maps in dimensions $n \geq 3$",
	\newblock \textit{Communications in Partial Differential Equations} \textbf{33} (2008) pp. 451 -- 477.

\bibitem{BIK}
	\newblock I. Bejenaru, A.D. Ionescu, and C. E. Kenig,
	\newblock ``Global existence and uniqueness of Schr{\"o}dinger maps in dimensions $d \geq 4$",
	\newblock \textit{Advances in Mathematics} \textbf{215} (2007) pp. 263 -- 291.
	
\bibitem{BIKT}
	\newblock I. Bejenaru, A.D. Ionescu, C. E. Kenig, and D. Tataru,
	\newblock ``Global Schr{\"o}dinger maps in dimensions $d \geq 2$ : small data in the critical Sobolev spaces",
	\newblock \textit{Annals of Mathematics} \textbf{173} (2) no. 3 (2011) pp. 1443 -- 1506.

\bibitem{B2}
	\newblock J. Bourgain,
	\newblock ``Refinements of {S}trichartz' inequality and applications to 2{D-NLS} with critical nonlinearity'',
	\newblock \textit {International Mathematical Research Notices}, \textbf{5} (1998) pp. 253 -- 283.
	
\bibitem{CGT}
	\newblock J. Colliander, M. Grillakis, and N. Tzirakis.
	\newblock ``Tensor products and correlation estimates with applications to nonlinear {S}chr{\"o}dinger equations'',
        \newblock \textit{Communications on Pure and Applied Mathematics}, \textbf{62} no. 7 (2009) pp. 920 -- 968.	
	
\bibitem{CKSTT2}
	\newblock J. Colliander, M. Keel, G. Staffilani, H. Takaoka, and T. Tao,
	\newblock ``Global existence and scattering for rough solutions of a nonlinear {S}chr{\"o}dinger equation on $\mathbf{R}^{3}$'',
        \newblock \textit{Communications on Pure and Applied Mathematics}, \textbf{21} (2004) pp. 987 - 1014.	

\bibitem{Dodson}
	\newblock B. Dodson,
	\newblock ``Bilinear Strichartz estimates for the Schr{\"o}dinger map problem",
	\newblock Arxiv eprints: 1210 : 5255 (2012).
	
\bibitem{DS}
	\newblock B. Dodson and P. Smith,
	\newblock ``A controlling norm for energy - critical Schr{\"o}dinger maps",
	\newblock \textit{Transactions of the American Mathematical Society} \textbf{367} (2015) pp. 7193 -- 7220.
	
\bibitem{KLPST}
	\newblock C. Kenig, T. Lamm, D. Pollack, G. Staffilani, and T. Toro,
	\newblock ``The Cauchy problem for Schr{\"o}dinger flows into K{\"a}hler manifolds",
	\newblock Preprint.		

\bibitem{Krieger}
	\newblock J. Krieger,
	\newblock ``Global regularity of wave maps from $\mathbf{R}^{3 + 1}$ to surfaces",
	\newblock \textit{Communications in Mathematical Physics} \textbf{238} (2003) pp. 333 -- 366.
	
\bibitem{Krieger1}
	\newblock J. Krieger,
	\newblock ``Global regularity of wave maps from $\mathbf{R}^{2 + 1}$ to $\mathbf{H}^{2}$. Small energy",
	\newblock \textit{Communications in Mathematical Physics} \textbf{250} (2004) pp. 507 -- 580.
	
\bibitem{NSU}
	\newblock A. Nahmod, A. Stefanov, and K. Uhlenbeck,
	\newblock ``On the well - posedness of the wave map problem in high dimensions",
	\newblock \textit{Comm. Anal. Geom.} \textbf{11} (2003) pp. 49 -- 83.	
	
\bibitem{PV}
	\newblock F. Planchon and L. Vega,
	\newblock ``Bilinear virial identities and applications''
	\newblock \textit{Annales Scientifiques de l'\'Ecole Normale Sup\'erieure Quatri\`eme S\'erie} \textbf{42} no. 2 (2009) pp. 261 - 290.	
	
\bibitem{ShSt}
	\newblock J. Shatah and M. Struwe,
	\newblock ``The Cauchy problem for wave maps",
	\newblock \textit{International Mathematics Research Notices} \textbf{2002} (2002) pp. 555 -- 571.
	
\bibitem{Smith1}
	\newblock P. Smith,
	\newblock ``Global regularity of critical Schr{\"o}dinger maps : subthreshold dispersed energy",
	\newblock Arxiv eprints : 1112.0251 (2011).
	
\bibitem{Smith2}
	\newblock P. Smith,
	\newblock ``Geometric renormalization below the ground state",
	\newblock \textit{International Mathematics Research Notices} \textbf{16} (2012) pp. 3800 -- 3844.	
	
\bibitem{Smith3}
	\newblock P. Smith,
	\newblock ``Conditional global regularity of Schr{\"o}dinger maps : subthreshold dispersed energy",
	\newblock \textit{Analysis of PDE} \textbf{6} (2013) pp. 601 -- 686.				
	
\bibitem{Taowm}
	\newblock T. Tao
	\newblock ``Global regularity of wave maps II. Small energy in two dimensions",
	\newblock \textit{Communications in Mathematical Physics} \textbf{224} (2001), pp. 443 -- 544.
	
\bibitem{TVZ}
	\newblock T. Tao, M. Visan, and X. Zhang.
	\newblock ``The nonlinear {S}chr\"odinger equation with combined power-type nonlinearities'',
	\newblock \textit{Communications in Partial Differential Equations}, \textbf{32} no. 7-9 (2007) pp. 1281--1343.	


\end{thebibliography}

\end{document}